\documentclass[a4paper, final, notitlepage, 10pt]{article}
\usepackage{latexsym,bm}
\usepackage{amsthm}
\usepackage{amsfonts}
\usepackage{amsmath}
\usepackage[all]{xy}
\usepackage[pdftex,final]{graphics}
\usepackage{fancyhdr}
\usepackage{verbatim}
\usepackage{geometry}
\usepackage{xcolor}
\usepackage{enumerate}
\usepackage{showkeys}

\newtheorem{thm}{Theorem}[section]

\newtheorem{lem}[thm]{Lemma}
\newtheorem{prop}[thm]{Proposition}
\theoremstyle{definition}

\newtheorem{rem}[thm]{Remark}

\numberwithin{equation}{section}

\newcommand{\thmref}[1]{Theorem~\ref{#1}}

\newcommand{\lemref}[1]{Lemma~\ref{#1}}
\newcommand{\propref}[1]{Proposition~\ref{#1}}

\newcommand{\remref}[1]{Remark~\ref{#1}}
\renewcommand{\O}{\mathcal{O}}
\renewcommand{\o}{\omega}
\newcommand{\E}{\mathcal{E}}

\renewcommand{\L}{\mathcal{L}}

\title{On canonically polarized Gorenstein $3$-folds \\satisfying the Noether equality}
\author{Yifan Chen~\qquad~Yong Hu}
\date{}
\begin{document}
\maketitle
\begin{abstract}
We study canonically polarized Gorenstein $3$-folds with at most terminal singularities and satisfying
$K_X^3=\frac 43p_g(X)-\frac {10}3$ and $p_g(X) \ge 7$.
We characterize the canonical maps of such $3$-folds,
describe a structure theorem for the locally factorial ones and completely classify the smooth ones.
New examples of canonically polarized smooth $3$-folds with $K_X^3=\frac 43p_g(X)-\frac {10}3$
and  $p_g(X) \ge 7$ are constructed. These examples are natural extensions of those constructed by M.~Kobayashi.
\end{abstract}
\section{Introduction}
Throughout the article, we work over the complex number field $\mathbb{C}$.

Let $S$ be a smooth minimal projective surface of general type.
We have the classical Noether inequality $K_{S}^2\geq 2 p_g(S)-4$ (cf.~\cite{noether}).
In \cite{horikawa} Horikawa classified surfaces on the Noether line.

Let $X$ be a projective $3$-fold of general type.
There have been many works dedicated to prove the $3$-dimensional version of the Noether inequality
(cf.~\cite{kobayashi}, \cite{mchen1}, \cite{mchen2}, \cite{CCZ} and \cite{chenchen}).
In \cite{chenchen},
the inequality $K_X^3 \ge \frac{4}{3}p_g(X)-\frac{10}{3}$ is proved
under the assumption that  $X$ is Gorenstein minimal.
This inequality is sharp according to Kobayashi's examples (cf.~\cite{kobayashi}).

In this article, we study the $3$-folds on the Noether line $K_X^3=\frac 43 p_g(X)-\frac{10}3$ and $p_g(X)\ge 7$.
We restrict our attention to the situation
where $X$ is a canonically polarized Gorenstein $3$-fold with at most terminal singularities,
and investigate such $3$-folds in two aspects: the classification and the construction of examples.
We first state the results of the construction of examples,
though it is the classification that leads to the discovery of the new ones.
\begin{thm}\label{thm:example}
Let $(e, a)$ be a pair of integers such that $a \ge e \ge 3$; or $1 \le e \le 2$, $a \ge e+1$; or $e=0$, $a \ge 2$.
Then there are  smooth $3$-folds $X$ with
$K_X^3=8a-4e-6$ and $p_g(X)=6a-3e-2$.
Moreover, the canonical divisor $K_X$ is ample and the canonical image of $X$ is the image of the embedding
 of the Hirzebruch surface $\Sigma_e$ into $\mathbb{P}^{p_g(X)-1}$ induced by the linear system $|s+(3a-e-2)l|$.
\end{thm}
We recall that $\Sigma_e$ is the projective bundle
$\mathbb{P}_{\mathbb{P}^1}(\O_{\mathbb{P}^1}\oplus \O_{\mathbb{P}^1}(-e))$ over $\mathbb{P}^1$.
Here $l$ stands for a fiber of the natural ruling $\Sigma_e \rightarrow \mathbb{P}^1$ and
$s$ stands for the section corresponding to the projection $\O_{\mathbb{P}^1} \oplus \O_{\mathbb{P}^1}(-e) \rightarrow \O_{\mathbb{P}^1}(-e)$.
It is clear that the $3$-folds $X$ are not isomorphic for different pairs $(e, a)$.
Our construction in Section~2 essentially follows the same method as \cite[(3.2)--(3.5)]{kobayashi}
and all these $3$-folds are finite double covers of certain $\mathbb{P}^1$-bundles over the Hirezbruch surfaces.
But we do construct more examples than \cite{kobayashi}.
Actually, the $3$-folds there correspond to the pairs $(e, a)$ with $a=e \ge 3$.
Also observe that $p_g(X) \ge 7$ for any $X$ in \thmref{thm:example}.
See \remref{rem:pg=4example} for examples of $3$-folds with $K_X^3=2$ and $p_g=4$, having exactly one singularity
and canonically fibred by curves.

Now we analyse the canonical maps of canonically polarized Gorenstein $3$-folds on the Noether line
and with $p_g(X) \ge 7$.
Assertions (a) and (b) in the following theorem characterize the base locus of $|K_X|$.
They are the key points in the classification.
\begin{thm}[cf.~(3.1)~Theorem in \cite{kobayashi}]\label{thm:baselocus}
Let $X$ be a Gorenstein $3$-fold with at most terminal singularities and satisfying
$K_X^3=\frac 43 p_g(X)-\frac{10}3$ and $p_g(X) \ge 7$.
Assume that $K_X$ is ample.
\begin{enumerate}[\upshape (a)]
\item The base locus of $|K_X|$ consists of a smooth rational curve $\overline{\Gamma}$,
      which is contained in the smooth locus of $X$.
\item Let $\pi \colon Y \rightarrow X$ be the blowup along $\overline{\Gamma}$
      and let $E_0$ be the exceptional divisor.
      Then $|\pi^*K_X-E_0|$ is base point free
      and it induces a fibration $\phi \colon Y \rightarrow \Sigma$,
      where $\Sigma$ is a surface in $\mathbb{P}^{p_g(X)-1}$ with $\deg \Sigma=p_g(X)-2$.
\item Let $C$ be the general fiber of $\phi$. Then $g(C)=2$ and $\pi^*K_X.C=1$.
\item The restriction $\phi|_{E_0} \colon E_0 \rightarrow \Sigma$ is birational.
\end{enumerate}
\end{thm}

To prove \thmref{thm:baselocus},
we need to reproduce the proofs of the Noether inequality
in \cite{mchen1}, \cite{mchen2} and \cite{chenchen}, and to get some geometric properties of the $3$-folds on the equality (see Section~3).
We remark that the \thmref{thm:baselocus}  no longer holds for the case where $p_g(X)=4$. See \remref{rem:pg=4example} for counterexamples and see
\remref{rem:pg=4baselocus} for further remarks.

To classify the $3$-folds in \thmref{thm:baselocus}, we shall determine
the surface $\Sigma$ explicitly and  analyse the fibration $\phi \colon Y \rightarrow \Sigma$.
From the classification of surfaces of minimal degree (for example, see~\cite[Exercises~IV.18~4)]{algebraicsurface} or \cite[p.~380, Corollary~2.19]{gtm52}), we know that $\Sigma$ is a Hirzebruch surface or a cone over a rational normal curve.
By showing the flatness of $\phi$ and by applying the semi-positivity theorem  to $\phi$ (cf.~\cite{viehweg}),
we obtain a structure theorem for $X$ when $X$ is locally factorial
and a complete classification when $X$ is smooth (See Section~4).

\begin{thm}\label{thm:classification}
Keep the same assumptions and notation in \thmref{thm:baselocus} and assume further that $X$ is locally factorial.
\begin{enumerate}[\upshape (a)]
\item  The surface $\Sigma$ is isomorphic to a Hirzebruch surface $\Sigma_e$ and $K_X^3=8a-4e-6$
       for a pair of integers $(e,a)$ such that $a \ge e \ge 3$; or $1 \le e \le 2$, $a \ge e+1$; or $e=0$, $a \ge 2$.
\item  The fibration $\phi \colon Y \rightarrow \Sigma$ is flat with irreducible fibers and
       $E_0$ is a section of $\phi$.
\end{enumerate}
Identify $\mathrm{Pic}(\Sigma)$ with $\mathrm{Pic}(\Sigma_e)=\mathbb{Z}s\oplus\mathbb{Z}l$ via $\Sigma \cong \Sigma_e$ in (a).
\begin{enumerate}[\upshape (c)]
\item  Then $\phi_\ast\O_Y(2E_0)=\O_{\Sigma}\oplus \O_\Sigma(-2s-2al)$ and $\phi_\ast\omega_{Y/\Sigma}=\O_{\Sigma}(3s+3al)\oplus \O_\Sigma(s+al)$,
       where $\omega_{Y/\Sigma}:=K_Y-\phi^*K_\Sigma$.
\end{enumerate}
Set $P:=\mathbb{P}_{\Sigma}(\O_{\Sigma}\oplus \O_\Sigma(-2s-2al))$.
       Denote by $\tau \colon P \rightarrow \Sigma$ the natural projection and by $E$ the section of $\tau$ corresponding to the projection
       $ \O_{\Sigma}\oplus \O_\Sigma(-2s-2al) \rightarrow \O_\Sigma(-2s-2al)$.
\begin{enumerate}[\upshape (d)]
\item   There is a finite double cover $\psi \colon Y \rightarrow P$ branched along $E+T$ such that $\phi=\tau\circ\psi$,
        where $T \in |5E+\tau^*(10s+10al)|$ and $T \cap E=\emptyset$.
\end{enumerate}
Moreover, if $X$ is smooth, then $X$ is one of the smooth $3$-folds constructed in Section~2.
\end{thm}

In \thmref{thm:classification}, we need the assumption that $X$ is locally factorial in order to apply intersection theory (cf.~\remref{rem:locallyfactorial}).
In generality, one would like to classify Gorenstein minimal $3$-folds on the Noether line.
However, this problem does not seem possible to resolve with the methods and the techniques of the present article.

\section{Examples}
Kobayashi constructed examples of $3$-folds on the Noether line in \cite{kobayashi}.
Using the same method there,
we construct more examples of canonically polarized smooth $3$-folds with $K_X^3=\frac 43p_g(X)-\frac{10}{3}$ and $p_g(X) \ge 7$.
See the diagram \eqref{eq:construction} for the  process.
We start with the Hirzebruch surfaces and certain $\mathbb{P}^1$-bundles over them
(cf.~\cite[p.~162, Proposition~7.11-7.12; p.~253, Exercise~8.4; Chapter~V, Section~2]{gtm52}).

\begin{align}\label{eq:construction}
\xymatrix{                                  & Y  \ar_{\pi}"2,1" \ar^{\psi}"1,3" \ar"2,2"     & P \ar_{\tau}"2,3" &   \\
\scriptstyle{\mathbb{P}(H^0(Y, \O_Y(mH))^*)}& X  \ar@{_{(}->}"2,1"             & \Sigma_e \ar"2,4" & \mathbb{P}^1
}
\end{align}

We denote by $\Sigma_e$ the Hirzebruch surface
$\mathbb{P}_{\mathbb{P}^1}(\O_{\mathbb{P}^1}\oplus \O_{\mathbb{P}^1}(-e))$ for $e \ge 0$.
Denote by $l$ the fiber of the natural ruling $\Sigma_e \rightarrow \mathbb{P}^1$ and
by $s$ the section corresponding to the projection of $\O_{\mathbb{P}^1} \oplus \O_{\mathbb{P}^1}(-e) \rightarrow \O_{\mathbb{P}^1}(-e)$.
Then
\begin{align}\label{eq:Hirzebruch}
s^2=-e,\ \O_{\Sigma_e}(s)=\O_{\Sigma_e}(1),\
\mathrm{Pic}(\Sigma_e)=\mathbb{Z}s\oplus\mathbb{Z}l\ \text{and}\ K_{\Sigma_e}=-2s-(e+2)l
\end{align}
The curve $s$ is called the negative section when $e >0$.

Let $\E$ be a locally free sheaf of rank $2$, sitting in the exact sequence
\begin{align}
    0 \rightarrow \O_{\Sigma_e} \rightarrow \E \rightarrow \O_{\Sigma_e}(-2s-2al) \rightarrow 0  \label{eq:rank2}
\end{align}
where $a$ is an integer, and let $P:=\mathbb{P}_{\Sigma_e}(\E)$.
Denote by $\tau \colon P \rightarrow \Sigma_e$ the natural projection and
by $E$ the section of $\tau$ corresponding to the morphism $\E \rightarrow \O_{\Sigma_e}(-2s-2al)$
in \eqref{eq:rank2}.
By abuse of notation,
we identify $\mathrm{Pic}(E)$ with $\mathbb{Z}s \oplus \mathbb{Z}l$ via $\tau|_{E} \colon E \cong \Sigma_e$.
\begin{lem}\label{lem:rank2}Keep the same notation as above. Then
\begin{enumerate}[\upshape (a)]
    \item $\O_P(E)=\O_P(1)$, $K_P=\tau^*(-4s-(2a+e+2)l)-2E$ and $\O_E(E)=\O_E(-2s-2al)$;
    \item the exact sequence \eqref{eq:rank2} splits if $a \ge e$;
    \item the linear system $|E+\tau^*(2s+2al)|$ is base point free if $a \ge e$.
\end{enumerate}
\end{lem}
\proof
The first two formulae of (a) are standard.
Since $K_E=-2s-(e+2)l$, the third one follows from the second one by the adjunction formula.

For (b), it suffices to show $H^1(\Sigma_e, \O_{\Sigma_e}(2s+2al))=0$.
According to \cite[p.~371, Lemma~2.4 and p.~162, Proposition~7.11]{gtm52},
$H^1(\Sigma_e, \O_{\Sigma_e}(2s+2al))=
H^1(\mathbb{P}^1, \O_{\mathbb{P}^1}(2a))\oplus H^1(\mathbb{P}^1, \O_\mathbb{P}^1(2a-e))\oplus
H^1(\mathbb{P}^1, \O_{\mathbb{P}^1}(2a-2e))$,
which is $0$ since $a \ge e$.

If $a \ge e$, then $|2s+2al|$ is base point free by \cite[p.~380, Corollary~2.8]{gtm52} and so is $\tau^*(2s+2al)$.
From the third formula in (a), we have $\O_P(E+\tau^*(2s+2al))|_E=\O_E$.
So the sequence
$0 \rightarrow \O_P(\tau^*(2s+2al)) \rightarrow \O_P(E+\tau^*(2s+2al)) \rightarrow \O_E \rightarrow 0$
is exact.
Since $R^i\tau_{\ast}\O_P=0$ for $i \ge 1$,
we have $H^1(P, \O_P(\tau^*(2s+2al))=H^1(\Sigma_e, \O_{\Sigma_e}(2s+2al)))=0$.
So the trace of $|E+\tau^*(2s+2al)|$ on $E$ is base point free.
This completes the proof of (c).\qed\smallskip

In the remaining of this section, we fix a pair $(e,a)$ as in the assumption of \thmref{thm:example}.
We use the theory of double covers to construct $3$-folds and prove \thmref{thm:example}.

Let $T$ be any smooth effective divisor in $|5E+\tau^*(10s+10al)|$.
Such $T$ exists by \lemref{lem:rank2}~(c) and the Bertini theorem.
Then $T \cap E = \emptyset$ since $\O_P(T)|_E=\O_E$ (see the proof of \lemref{lem:rank2}~(c)).
We also have
\begin{align}
E+T \sim 2\L,\ \text{where}\ \L=3E+\tau^*(5s+5al) \label{eq:coverdata}
\end{align}
Therefore there is a finite smooth double cover $\psi \colon Y \rightarrow P$ branched along $E+T$
such that $\psi_{\ast}\O_Y = \O_P \oplus \L^{\vee}$.
Set $E_0:=\psi^{-1}(E)$. Then $\psi^*E=2E_0$.
We identify $\mathrm{Pic}(E_0)$ with $\mathbb{Z}s \oplus \mathbb{Z}l$ via
$\tau|_E \circ \psi|_{E_0} \colon E_0 \cong E \cong \Sigma_e$.

\begin{prop}\label{prop:coverinvariant}
Let $H:=K_Y-E_0$. Then
\begin{enumerate}[\upshape (a)]
\item $K_Y=\psi^*\tau^*(s+(3a-e-2)l)+2E_0$ and $p_g(Y)=6a-3e-2$;
\item $H$ is nef and $H^3=8a-4e-6$;
\item for an irreducible curve $C$ in $Y$,
      $HC=0$ holds if and only if $C$ is a fiber of the ruling of $E_0$ induced by $|l|$;
\item $3H-K_Y$ is ample.
\end{enumerate}
\end{prop}
\proof
For (a), we have $K_Y=\psi^*(K_P+\L)$ and $p_g(Y)=p_g(P)+h^0(P, \O_P(K_P+\L))=h^0(P, \O_P(K_P+\L))$ since $p_g(P)=0$.
Note that $K_P+\L=E+\tau^*(s+(3a-e-2)l)$ by the second formula of \lemref{lem:rank2}~(a) and \eqref{eq:coverdata}.
The formula for $K_Y$ follows.
By \lemref{lem:rank2}~(a)-(b),
we have
$\tau_{\ast}\O_P(K_P+\L)=\E \otimes \O_{\Sigma_e}(s+(3a-e-2)l)=\O_{\Sigma_e}(s+(3a-e-2)l)\oplus \O_{\Sigma_e}(-s+(a-e-2)l)$.
As in the proof of \lemref{lem:rank2}~(b),
a similar calculation yields
$h^0(\Sigma_e, \O_{\Sigma_e}(s+(3a-e-2)l))=6a-3e-2$ and $h^0(\Sigma_e, \O_{\Sigma_e}(-s+(a-e-2)l))=0$.
Therefore $p_g(Y)=h^0(P, \tau_{\ast}\O_P(K_P+\L))=6a-3e-2$.

Let $M:=\psi^*\tau^*(s+(3a-e-2)l)$.
Then $H=M+E_0$, $K_Y=M+2E_0$ and $M|_{E_0}=s+(3a-e-2)l$.
So $E_0|_{E_0}=\frac 13(K_{E_0}-M|_{E_0})=-s-al$ by the adjunction formula and then $H|_{E_0}=(2a-e-2)l$.
Therefore $H^2E_0=0$, $HME_0=H|_{E_0}.M|_{E_0}=2a-e-2$ and $M^2E_0=(M|_{E_0})^2=6a-3e-4$.
It follows that $H^3=8a-4e-6$.

Now assume $H.C \le 0$ for an irreducible curve $C$.
If $\tau\psi(C)$ is a point, then $H.C=\frac 12 \psi^*E.C >0$ since $E$ is a section of $\tau$.
So $\tau\psi(C)$ is a curve.
Since $s+(3a-e-2)l$ is very ample by \cite[p.~380, Corollary~2.18]{gtm52},
$MC>0$ and thus $E_0C=(H-M)C<0$.
Therefore $C$ is contained in $E_0$.
Since $H|_{E_0}=(2a-e-2)l$ and $2a-e-2>0$,
we conclude that $H.C=0$ and $C$ is a fiber of the ruling of $E_0$.
This completes the proof of (b) and (c).

Because both $H$ and $M$ are nef and $H^3 >0$,  $3H-K_Y=H+M$ is nef and big.
Since $M|_{E_0}$ is very ample, it follows by (c) that $(3H-K_Y)C>0$ for any curve $C$.
Assume $(3H-K_Y)^2S=(H+M)^2S=0$ for some irreducible surface $S$.
Then $H^2S=HMS=M^2S=0$ and thus $ME_0S=M(H-M)S=0$.
Because $M|_{E_0}$ is very ample, we conclude $S \cap E_0= \emptyset$.
Then the equality $HMS=0$ yields a contradiction to (c), since $|M|$ is base point free.
We have shown $(3H-K_Y)^2S>0$ for any irreducible surface $S$.
Therefore $3H-K_Y$ is ample.
\qed\smallskip

By \propref{prop:coverinvariant}~(c)-(d),
the base-point-free theorem \cite[Theorem~3-1-1]{KMM} implies that the image $X$ of the morphism
$\pi \colon Y \rightarrow \mathbb{P}(H^0(Y, \O_Y(mH))^*)$ for $m>>0$
can be identified with the blowdown along the ruling of $E_0$.
Therefore $X$ is a $3$-fold birational to $Y$.

We now show that $X$ satisfies the properties in \thmref{thm:example}.
Observe that $(3H-K_Y)+tK_Y|_{E_0}=3l+(t-1)(-s+(a-e-2)l)$ for any $t \in \mathbb{R}$.
Hence $\mathrm{max}\{t\in \mathbb{R}| 3H-K_Y+tK_Y\ \text{is nef}\}=1$.
Since  $\O_{E_0}(E_0)\cong \O_{E_0}(-s-al)$, $X$ is smooth by \cite[Theorem~(3.3)~(3.3.1)]{mori}.
Also $K_X$ is ample by the base-point-free theorem.
From the construction,  we have $\pi^*K_X \cong H$ and thus $K_X^3=H^3=8a-4e-6$.
Since $X$ is birational to $Y$, $p_g(X)=p_g(Y)=6a-3e-2$.
We have seen that $H=E_0+\psi^*\tau^*(s+(3a-e-2)l)$ in the proof of \propref{prop:coverinvariant}.
Because the movable part of $|\pi^*K_X|=|H|$ is $\psi^*\tau^*(s+(3a-e-2)l)$,
we see that the base locus $|K_X|$ consists of the smooth rational curve $\pi(E_0)$
and the canonical image of $X$ is the image of the embedding of $\Sigma_e$ into $\mathbb{P}^{p_g(X)-1}$
 induced by $|s+(3a-e-2)l|$.

\begin{rem}\label{rem:pg=4example}
If $(e, a)=(0, 1)$,
the construction of $Y$ still works and \propref{prop:coverinvariant} still holds except the assertion~(c).
Indeed, $H|_{E_0}=0$ in this case and $|mH|$ ($m >>0$) contracts exactly the whole divisor $E_0$.
The $3$-fold $X$ obtained by the base-point-free theorem has ample canonical divisor
and it is still on the Noether line.
Indeed, $K_X^3=2$ and $p_g(X)=4$.
But $X$ is no longer smooth.
It has exactly one singularity (cf.~\cite[Theorem~(3.3)~(3.3.3)]{mori}).
Moreover, $|K_X|$ has this singularity as base locus and the canonical image of $X$ is a smooth quadric.
\end{rem}

\begin{rem}\label{rem:gorensteinexample}
In the construction above, if we allow $T$ ($\in |5E+\tau^*(10s+10al)|$) to be singular, then so is $X$.
However, it is difficult to explicitly construct $T$ with  prescribed singularities in  the linear system $|5E+\tau^*(10s+10al)|$.
\end{rem}
\section{Base locus of the canonical linear systems}
We prove \thmref{thm:baselocus} in this section.
Throughout this section,
we denote by $X$ a canonically polarized Gorenstein $3$-fold with at most terminal singularities.
In order to prove \thmref{thm:baselocus}, we need to reproduce the proofs of the Noether inequality
in \cite{mchen1}, \cite{mchen2} and \cite{chenchen}.

\subsection{Setting}
This subsection is devoted to study the canonical map of $X$.
Write
\begin{align}\label{eq:movablefixed}
|K_X|=|\overline{M}|+\overline{Z}
\end{align}
where $|\overline{M}|$ is the movable part of $|K_X|$ and $\overline{Z}$ is the fixed part of $|K_X|$.

We shall resolve the base locus of $|\overline{M}|$ in two steps.
For a linear system $\Upsilon$, we denote by $\mathrm{Bs}\Upsilon$ the base locus of $\Upsilon$.
Roughly speaking,
the first step is to resolve the subset $\mathrm{Bs}|\overline{M}| \cap \mathrm{Sing}(X)$.
\begin{prop}[{cf.~\cite[Section~2]{chenchen}}]\label{prop:fgresolution}
There is a birational morphism $\alpha \colon X_0 \rightarrow X$
satisfying the following properties.
\begin{enumerate}[\upshape (a)]
\item The morphism $\alpha$ is a composition of successive divisorial contractions to points and $X_0$ is a Gorenstein $3$-fold with at most terminal singularities.
\item Denote by $|M_0|$  the movable part of $|\alpha^*\overline{M}|$.
      Then $\mathrm{Bs}|M_0|\cap \mathrm{Sing}(X_0)=\emptyset$.
\item The following formulae
      \begin{align}
      K_{X_0}=\alpha^*K_X+\sum_{t=1}^{m}c_tD_t,\ \alpha^*(\overline{M})=M_0+\sum_{t=1}^{m}d_tD_t,\ \alpha^*(\overline{Z})=Z_0+\sum_{t=1}^{m}e_tD_t \label{eq:fGadjunction}
      \end{align}
      hold, where
      \begin{enumerate}[\upshape (i)]
      \item $Z_0$ is the strict transform of $\overline{Z}$,
         \item $D_t$ is a $\mathbb{Q}$-Cartier prime divisor such that $\alpha(D_t)$ is a point for $1 \le t \le m$, and
         \item $c_t$ and $d_t+e_t$ are integers such that $0 <c_t \le d_t+e_t$ for $1\le t \le m$.
      \end{enumerate}
\item Set $D:=\cup_{t=1}^m D_t$. Then $D$ is a disjoint union of fibers of $\alpha$.
\end{enumerate}
\end{prop}
\proof The birational morphism $\alpha$ is  constructed in \cite[p.~4--p.~5]{chenchen},
using explicit resolutions of terminal singularities (see \cite{jachen} and \cite[Definition~2.2]{chenchen}).
Then (a) and (b) follow from the construction.
Since both $X_0$ and $X$ are Gorenstein, $c_t$ and $d_t+e_t$ are integers.
The inequality $c_t\le d_t+e_t$ follows by \cite[Corollary~2.4]{chenchen}.
Assertion (d) follows from the fact that $\alpha$ is a composition of successive divisorial contractions to points.\qed

We fix a birational morphism $\alpha \colon X_0 \rightarrow X$ as in \propref{prop:fgresolution}.
We may assume that the number of divisorial contractions in the construction of $\alpha$ is minimal.
The second step is to resolve the base locus $|M_0|$ without changing the singularities of $X_0$.
This is possible by \propref{prop:fgresolution}~(b) and by Hironaka's Theorem (cf.~\cite{hironaka}).
\begin{prop}[{cf.~\cite[Lemma~4.2]{mchen1}}]\label{prop:blowupresolution}
There are successive blowups
\begin{align*}
 \beta \colon Y=X_{n+1}\stackrel{\pi_n}\rightarrow X_{n}\rightarrow \cdots \rightarrow X_{i+1} \stackrel{\pi_{i}}\rightarrow X_{i}\rightarrow \cdots \rightarrow X_1\stackrel{\pi_0}\rightarrow X_0
\end{align*}
such that $\pi_i$ is a blowup along a smooth irreducible center $W_i$,
$W_i$ is contained in the base locus of the movable part of $|(\pi_0\circ \pi_1 \circ \cdots \circ \pi_{i-1})^*M_0|$ and $W_i \cap \mathrm{Sing}(X_i)=\emptyset$.
Moreover, the morphism $\beta=\pi_n\circ \cdots \circ \pi_0$ satisfies the following properties.
\begin{enumerate}[\upshape (a)]
\item Denote by $|M|$ the movable part of $|\beta^*M_0|$. Then $|M|$ is base point free.
\item The following formulae
      \begin{align}\label{eq:blowupadjunction}
          K_{Y}=\beta^*K_{X_0}+\sum_{i=0}^{n}a_i E_i,\ \beta^*M_0=M+\sum_{i=0}^{n}b_i E_i
      \end{align}
      hold, where $E_i$ is the strict transform of the exceptional divisor of $\pi_i$ for $0 \le i \le n$, and
      $a_i$ and $b_i$ are positive integers such that $a_i \le 2b_i\ \text{for}\ 0 \le i \le n$.
\item If $a_k=b_k=1$ for some $k$ such that $0 \le k \le n$, then $W_k$ is a curve of $X_k$ and the general member of $|M_0|$ is smooth
      at a general point of the curve $\pi_{0} \circ \cdots \pi_{k-2} \circ \pi_{k-1}(W_k)$.
\item If $a_k=2b_k$ for some $k$ such that $0 \le k \le n$, then $W_k$ is a point of $X_k$ and
      the general member of $|M_0|$ is smooth at the point $\pi_{0}\circ\cdots \pi_{k-2}\circ\pi_{k-1}(W_k)$.
\end{enumerate}
\end{prop}
\proof
The construction of the blowups $\pi_i$ and  (a) follow by \propref{prop:fgresolution}~(b) and by Hironaka's Theorem (cf.~\cite{hironaka}).
We remark that the assertion $a_i \leq 2b_i$ in (b) is exactly  \cite[Lemma 4.2]{mchen1}.
What really involved here is the assertions~(c) and (d).

We introduce some notation.
For any $0 \le j<i \le n$,
let $E_{j}^{(i)}$ $(\subset X_i)$ be the strict transform of the exceptional divisor of $\pi_j$
and let $E_j^{(n+1)}:=E_j$.
According to the definition of $a_i$ and $b_i$, for any $0 \le i \le n$,
we have
\begin{align}
&K_{X_{i+1}}=(\pi_0\circ \pi_1 \circ \cdots \circ \pi_i)^*K_{X_0}+\sum_{j=0}^{i}a_j E_j^{(i+1)}, \nonumber\\
&(\pi_0\circ \pi_1 \circ \cdots \circ \pi_i)^*(M_0)=M_{i+1}+\sum_{j=0}^{i}b_j E_{j}^{(i+1)}, \label{eq:km}
\end{align}
where $|M_{i+1}|$ is the movable part of $|K_{X_{i+1}}|$.

Considering the single blowup $\pi_k \colon X_{k+1} \rightarrow X_k$ for $1 \le k \le n$,
we have for $j<k$,
\begin{align*}
&\pi_{k}^*E_j^{(k)}=E_j^{(k+1)}+r_jE_k^{(k+1)},\
\text{where $r_j \in \mathbb{Z}$ and $r_j \ge 0$, and $r_j >0$ if and only if $W_k \subseteq E_{j}^{(k)}$},\\
&K_{X_{k+1}}=\pi_k^*K_{X_k}+a_k'E_k^{(k+1)},\
\text{where $a_k'=1$ if $W_k$ is a curve or $a_k'=2$ if $W_k$ is a point},\\
&\pi_k^*M_k=M_{k+1}+b_k'E_k^{(k+1)},\ \text{where $b_k' \ge 1$ since $W_k$ is contained in the base locus of $|M_k|$.}
\end{align*}
Comparing these formulae with \eqref{eq:km} when $i=k$ and $i=k-1$, we obtain
\begin{align}
a_k=a_k'+\sum_{j=0}^{k-1}r_ja_j,\ b_k=b_k'+\sum_{j=0}^{k-1}r_jb_j\ \text{for}\ 1 \le k \le n \label{eq:aa'bb'}
\end{align}

Using \eqref{eq:aa'bb'}, we easily conclude that $a_i \le 2b_i$ for $0 \le i \le n$ by induction on $i$ and (b) follows.

For (c), the case $k=0$ is trivial.
If $k>0$, by \eqref{eq:aa'bb'}, we have $a_k'=1$, $b_k'=1$ and $r_j=0$ for any $j<k$.
Therefore $W_k$ is a smooth curve of $X_k$,
the general member of $|M_k|$ is smooth at a general point of $W_k$
and  $W_k \not \subset E_j^{(k)}$ for $j<k$.
Assertion (c) follows.

We proceed to prove (d) by induction on $k$.
Assume $a_0=2b_0$.
Note that $a_0=1$ if $W_0$ is a curve or $a_0=2$ if $W_0$ is a point.
Therefore $a_0=2$ and $b_0=1$.
It follows that $W_0$ is a point of $X$ and the general member of $|M_0|$ is smooth at $W_0$.

Now assume $a_k=2b_k$ for some $k>0$.
Since $a_j \le 2b_j$ for any $j<k$, by \eqref{eq:aa'bb'}, we have
$a_k'=2$, $b_k'=1$ and $a_j=2b_j$ for those $j$ such that $j<k$ and $r_j>0$.


Therefore $W_k$ is a smooth point of $X_k$ and the general member of $M_k$ is smooth at $W_k$.
If $r_j=0$ for any $j<k$,
then the point $W_k \not \in E_j^{(k)}$ for any $j<k$
and thus the second statement of (b) clearly holds.
If $r_j >0$ for some $j<k$,
then $a_j=2b_j$ and $\pi_{0}\circ\cdots \pi_{k-2}\circ\pi_{k-1}(W_k)=\pi_{0}\circ\cdots \pi_{j-2}\circ\pi_{j-1}(W_j)$.
Assertion (d) follows by induction.\qed

From now on, we fix a birational morphism $\beta$ as in \propref{prop:blowupresolution} such that the number $n+1$ of blowups is minimal.
Denote by $\phi_{K_X}$ the canonical map of $X$ and by $\Sigma$ the image of $\phi_{K_X}$.
Let $\phi$ be the morphism induced by the linear system $|M|$.
Then $\phi = \phi_{K_X} \circ \pi$, where $\pi=\alpha \circ \beta$.
Let $Y\stackrel{f}\rightarrow B \stackrel{\delta} \rightarrow \Sigma$ be the Stein factorization of $\phi$.
We have the following commutative diagram:
\begin{align*}
\xymatrix{
 & Y \ar[d]_{\pi} \ar[r]^{f} \ar[dr]^{\phi} \ar[dl]_{\beta}
                & B \ar[d]^{\delta} \\
X_0 \ar[r]^{\alpha} & X  \ar@{.>}[r]_{\phi_{K_X}}
                & \Sigma          }
\end{align*}
Note that $B$ is normal.
We have the following known results:
if $\dim B=3$, then $K_X^3\geq 2p_g(X)-6$ (cf.~\cite[Main Theorem]{kobayashi});
if $\dim B=2$, then $K_X^3\geq \left\lceil \frac{2}{3}(g(C)-1)\right\rceil(p_g(X)-2)$
where $g(C)$ is the genus of a general fiber $C$ of $f$ (cf.~\cite[Theorem 4.1~(ii)]{mchen1});
if $\dim B=1$, then $K_X^3\geq \frac{7}{5}p_g(X)-2$
(cf.~\cite[Theorem 4.1~(iii)]{mchen1} and \cite[Theorem 4.1]{CCZ}).

\subsection{Some basic results}
From now on,
we assume further that $X$ is on the Noether line $K_X^3=\frac 43 p_g(X)-\frac{10}{3}$ and $p_g(X) \ge 7$.

\begin{lem}[{cf.~\cite[Section~3]{kobayashi}}]\label{lem:fibration}Keep the same notation as above.
\begin{enumerate}[\upshape (a)]
    \item The morphism $f$ is a fibration over the normal surface $B$.
    \item Let $C$ be a general fiber of $f$.
          Then $C$ is a smooth curve of genus $2$ with $\pi^*K_X.C=1$.
    \item Let $d_{\Sigma}$ be the degree of $\Sigma$ in $\mathbb{P}^{p_g(X)-1}$.
          Then $d_{\Sigma} \ge p_g(X)-2$ and $M^2 \equiv d_{\Sigma}C$.
    \item The morphism $\delta$ is birational.
\end{enumerate}
\end{lem}
Here the symbol $\equiv$ in (c) stands for numerical equivalence.
\proof
From the discussion at the end of the last subsection,
we see that $B$ is a normal surface and $g(C)=2$.
In particular, we have $M^2 \equiv d_{\Sigma}\cdot \deg\delta \cdot C$.

Because $\Sigma$ is non-degenerate,  we have $d_{\Sigma} \ge p_g(X)-2$.
Since both $\pi^*K_X$ and $M$ are nef, we conclude that
$$K_X^3 \ge \pi^*K_XM^2 =d_{\Sigma}\cdot \deg\delta \cdot \pi^*K_X.C \ge (p_g(X)-2)\deg\delta \cdot \pi^*K_X.C.$$
If $\pi^*K_X.C\ge 2$ or $\deg\delta \ge 2$,  then $K_X^3 \ge 2p_g(X)-4$,
a contradiction to the Noether equality.
Therefore $\pi^*K_X.C=1$, $\deg\delta=1$ and thus $M^2\equiv d_{\Sigma}C$.\qed

In what follows, we shall figure out the geometric information hidden in the Noether equality $K_X^3=\frac 43p_g(X)-\frac {10}{3}$.
For this purpose, we use the techniques in the proofs of \cite[Theorem~4.3]{mchen1} and \cite[Theorem~3.1]{chenchen}.
Recall from \eqref{eq:movablefixed}, \eqref{eq:fGadjunction} and \eqref{eq:blowupadjunction} that
\begin{align}\label{eq:fulladjunction}
K_Y=\pi^*K_X+(\sum_{t=1}^{m}c_t\beta^*D_t+\sum_{i=0}^{n}a_iE_i),&&
\pi^*K_X=M+(\sum_{t=1}^{m}(d_t+e_t)\beta^*D_t+\sum_{i=0}^{n}b_iE_i+\beta^*{Z_0})
\end{align}
In particular, we have
$K_X^3=(\pi^*K_X)^3=(\pi^*K_X)^2M+K_X^2\overline{Z}$.
We aim to bound $(\pi^*K_X)^2M$ from below.

For this purpose, by Bertini's theorem, we choose a general member $S$ of $|M|$ such that $S$ is smooth and consider the fibration $f|_S$.
By abuse of notation, we still denote by $C$ the general fiber of $f|_S$.
We remark that $S$ is of general type since so is $X$.
Also the divisors $\beta^*D_t|_S$ and $E_i|_S$ are effective for $1 \le t \le m$ and $0 \le i \le n$.
In the following lemma, we decompose the effective divisors $(\sum_{t=1}^{m}c_t\beta^*D_t+\sum_{i=0}^{n}a_iE_i)|_S$
and $(\sum_{t=1}^{m}(d_t+e_t)\beta^*D_t+\sum_{i=0}^{n}b_iE_i+\beta^*{Z_0})|_S$ into the horizontal parts and the vertical parts with respect to the fibration $f|_S$.
\begin{lem}[{See the proof of ~\cite[Theorem~4.3]{mchen1}}]\label{lem:section}
We have
\begin{align}
(\sum_{t=1}^{m}c_t\beta^*D_t+\sum_{i=0}^{n}a_iE_i)|_S&=\Gamma+D_V+E_V,\nonumber\\
(\sum_{t=1}^{m}(d_t+e_t)\beta^*D_t+\sum_{i=0}^{n}b_iE_i+\beta^*{Z_0})|_S&=\Gamma+D_V'+E_V'\label{eq:vertical}
\end{align}
where $\Gamma$ is an irreducible reduced curve  with $\Gamma C=1$,
while $E_V, D_V, E_V'$ and $D_V'$ are effective divisors contained in the fibers of $f|_S$.
Moreover, exactly one of the following two cases occurs:
\begin{enumerate}[\upshape (i)]
    \item $\Gamma \subseteq \mathrm{Supp}(\sum_{t=1}^{m}\beta^*D_t|_S)$ ; in this case,
          $D_V=\sum_{t=1}^{m}c_t\beta^*D_t|_S-\Gamma$, $E_V=\sum_{i=0}^{n}a_iE_i|_S$,
          $D_V'=\sum_{t=1}^{m}(d_t+e_t)\beta^*D_t|_S-\Gamma$ and $E_V'=\sum_{i=0}^{n}b_iE_i|_S+\beta^*{Z_0}|_S$.
    \item $\Gamma \le E_l|_S$ for a unique integer $l$; in this case, $a_l=b_l=1$,
          $D_V=\sum_{t=1}^{m}c_t\beta^*D_t|_S$, $E_V=\sum_{i=0, i\not=l}^{n}a_iE_i|_S+(E_l|_S-\Gamma)$, $D_V'=\sum_{t=1}^{m}(d_t+e_t)\beta^*D_t|_S$
           and $E_V'=(\sum_{i=0, i\not=l}^{n}b_iE_i|_S+\beta^*{Z_0}|_S)+(E_l|_S-\Gamma)$.
\end{enumerate}
In either case, we have $D_V\le D_V'$ and $E_V \le 2E_V'$.
\end{lem}
\proof
Since $g(C)=2$, $K_Y\cdot C=2$ by the adjunction formula.
According to \lemref{lem:fibration}~(b) and \eqref{eq:fulladjunction},
we have
\begin{align*}
(\sum_{t=1}^{m}c_t\beta^*D_t+\sum_{i=0}^{n}a_iE_i)|_S\cdot C=1\ \text{and}\
(\sum_{t=1}^{m}(d_t+e_t)\beta^*D_t+\sum_{i=0}^{n}b_iE_i+\beta^*{Z_0})|_S\cdot C=1.
\end{align*}
We conclude that the horizontal part of $(\sum_{t=1}^{m}c_t\beta^*D_t+\sum_{i=0}^{n}a_iE_i)|_S$ consists of an irreducible reduced curve $\Gamma$, which is also the horizontal part of $(\sum_{t=1}^{m}(d_t+e_t)\beta^*D_t+\sum_{i=0}^{n}b_iE_i+\beta^*{Z_0})|_S$ and $\Gamma\cdot C=1$.

From \propref{prop:fgresolution} and \propref{prop:blowupresolution},
we see that $c_t, d_t+e_t, a_i$ and $b_i$ are positive integers such that $c_t \le d_t+e_t$ and $a_i \le 2b_i$.
Moreover, $D_t$ is $\mathbb{Q}$-Cartier and $E_i$ is Cartier.

If $E_lC>0$ for some $l$, then $E_lC=1$ since $E_i$ is Cartier, and then $\Gamma \le E_l|_S$.
We see that the case (ii) occurs.
If $E_iC=0$ for all $i$, then $(\sum_{t=1}^{m}c_t\beta^*D_t)C=1$ and $\Gamma \le \sum_{t=1}^{m}\beta^*D_t$,
and the case (i) occurs.

The last statement follows easily from (i), (ii) and the inequalities $c_t \le d_t+e_t$ and $a_i \le 2b_i$.\qed

\begin{rem}\label{rem:Esupport}
We remark that if $\overline{Z}=0$, then $\mathrm{Supp}(E_V)=\mathrm{Supp}(E_V')$ in either case.
This is indeed true by the following proposition.
Also one sees that $\mathrm{Supp}(D_V)=\mathrm{Supp}(D_V')$ in the case (ii).
We also remark that the case (ii) occurs, again by the following proposition.
\end{rem}

\begin{prop}[{See the proofs of \cite[Theorem~4.3]{mchen1} and \cite[Theorem~3.1]{chenchen}}]\label{prop:equalities}
In the above setting, we have
\begin{enumerate}[\upshape (a)]
    \item $\overline{Z}=0$;
    \item $\mathrm{Supp}(E_V)=\mathrm{Supp}(E_V')$ and $D_V\Gamma=D_V'\Gamma=(2E_V'-E_V)\Gamma=0$;
    \item $\delta \colon B \rightarrow \Sigma$ is an isomorphism
          and $\Sigma$ is a normal rational surface with $d_{\Sigma}=p_g(X)-2$;
    \item $\Gamma$ is a smooth rational curve and $\pi^*K_X|_S. \Gamma=\frac 13p_g(X)-\frac 43$;
    \item $\pi^*K_X|_S.E_V=\pi^*K_X|_S.E_V'=\pi^*K_X|_S\cdot D_V=\pi^*K_X|_S\cdot D_V'=0$.
    \item $\Gamma \le E_l|_S$ for a unique integer $l$, and for this integer $l$,
          $a_l=b_l=1$, $E_lC=1$ and $E_iC=0$ for $i\not=l$.
\end{enumerate}
\end{prop}
\proof
The last statement of \lemref{lem:section} yields that
$D_V'-D_V \ge 0$, $2E_V'-E_V\ge0$ and both divisors are contained in fibers of of $f|_S$.
Because $\Gamma$ is a section of $f|_S$,
\begin{align}\label{eq:gammaEV}
\Gamma(D_V'-D_V)\geq 0,\ \Gamma(2E'_V-E_V)\geq 0
\end{align}

The adjunction formula yields
\begin{align}\label{eq:gammagenus}
K_S\Gamma+\Gamma^{2}=2p_a(\Gamma)-2 \ge -2
\end{align}
Note that $K_Y|_S=\pi^*K_X|_S+\Gamma+D_V+E_V$ and $\pi^*K_X|_S=M|_S+\Gamma+D_V'+E_V'$ by \eqref{eq:fulladjunction} and \lemref{lem:section}.
By \eqref{eq:gammaEV}, one has
\begin{align}
-2 \le (K_S+\Gamma)\Gamma&=(K_Y+M)|_S\Gamma+\Gamma^2 \nonumber\\
                  &=(\pi^*K_X|_S+M|_S+D_V+E_V+2\Gamma)\Gamma \nonumber\\
                  &\le (\pi^*K_X|_S+M|_S+2D_V'+2E_V'+2\Gamma)\Gamma \label{eq:gammaEVeq}\\
                  &=(3\pi^*K_X|_S-M|_S)\Gamma \nonumber\\
                  &=3\pi^*K_X|_S.\Gamma-d_{\Sigma}.\nonumber
\end{align}
The last equality holds by $M^2 \equiv  d_{\Sigma}C$ (see \lemref{lem:fibration}~(c)) and $\Gamma.C=1$.
Also $d_{\Sigma} \ge p_g(X)-2$, we obtain
\begin{align}\label{eq:canonicalgamma}
\pi^*K_X|_S.\Gamma \ge \frac 13 (d_{\Sigma}-2) \ge \frac 13 (p_g(X)-4)
\end{align}

Finally, we have
\begin{align}
K_X^3=(\pi^*K_X)^3&=(\pi^*K_X|_S)^2+K_X^2\overline{Z} \nonumber\\
                  &=\pi^*K_X|_S.M|_S+\pi^*K_X|_S.\Gamma+\pi^*K_X|_S.D_V'+\pi^*K_X|_S.E_V'+K_X^2\overline{Z} \nonumber\\
                  &\ge d_{\Sigma}+(\frac 13 p_g(X)-\frac 43) +\pi^*K_X|_S.D_V'+\pi^*K_X|_S.E_V'+K_X^2\overline{Z} \nonumber\\
                  &\ge \frac 43p_g(X)-\frac{10}3 \label{eq:estimate}
\end{align}
By assumption, we see that all the equalities in the inequalities \eqref{eq:gammaEV}-\eqref{eq:estimate} hold.

Since $K_X$ is ample, (a) follows by the equality $K_X^2\overline{Z}=0$ in \eqref{eq:estimate}. Therefore we have $Z_0=0$ and $e_t=0$ for all $1\leq t\leq m$.
Then (b) follows by \remref{rem:Esupport}, the equalities $\Gamma(D_V'-D_V)=\Gamma(2E'_V-E_V)=0$ in \eqref{eq:gammaEV} and the one $D_V\Gamma=2D'_V\Gamma$ in \eqref{eq:gammaEVeq}.

For (c), we conclude $d_{\Sigma}=p_g(X)-2$ from \eqref{eq:canonicalgamma}.
This implies that $\Sigma$ is of minimal degree.
Since  $\Sigma$ is non-degenerate, $\Sigma$ is a normal rational surface
(cf.~\cite[Exercises IV.~18.~4]{algebraicsurface}).
Because $B$ is also normal and $\delta$ is finite by Stein factorization,
$\delta$ is an isomorphism by \lemref{lem:fibration}~(d).

Assertion (d) follows by the equalities in \eqref{eq:gammagenus} and \eqref{eq:canonicalgamma}.
For (e), the equality in \eqref{eq:estimate} implies $\pi^*K_X|_S\cdot D_V'=\pi^*K_X|_S \cdot E_V'=0$.
Since $\pi^*K_X|_S$ is nef,  $D_V \le D_V'$  and
$E_V \le 2E_V'$  by \lemref{lem:section},
we conclude $\pi^*K_X|_S\cdot D_V=\pi^*K_X|_S \cdot E_V=0$.

Since $p_g(X)\geq 7$, $\pi(\Gamma)$ is a curve by (d).
Because $\alpha(D_t)$ is a point for $1 \le t \le m$ (see \propref{prop:fgresolution}) and $\pi=\alpha \circ \beta$,
we conclude that the case (ii) of \lemref{lem:section} occurs.
It remains to show $E_lC=1$ and $E_iC=0$ for $i\not= l$ in (f).
By \lemref{lem:section}, $E_VC=0$ and $E_V=(\sum_{i=0, i\not=l}^{n}a_iE_i)|_S+(E_l|_S-\Gamma)$.
Because $a_i$ is positive and $\Gamma C=1$, we obtain the required equalities.
\qed

\begin{rem}\label{rem:pg=4baselocus}
In the discussion above, we use the assumption $p_g(X) \ge 7$ only in the proof of \lemref{lem:fibration}, to show that  the canonical image of $X$ is a surface,
and in the proof of \propref{prop:equalities}~(f), to show that $\pi(\Gamma)$ is curve.
The fact that $\pi(\Gamma)$ is a curve plays an important role in the proof of \thmref{thm:baselocus}.
For example,
it is the key point to show that we actually do not need the first step to resolve $\mathrm{Bs}|\overline{M}|$ (see \lemref{lem:alpha} and \propref{prop:fgresolution}).

It is proved in \cite{kobayashi} that $K_X^3 \ge 2p_g(X)-6$
if the $3$-fold $X$ has at most canonical singularities and $\dim \mathrm{Im}\phi_{K_X}=3$.
Note that the lines $K_X^3=2p_g(X)-6$ and $K_X^3=\frac 43p_g(X)-\frac{10}{3}$
intersect at a point $(K_X^3,p_g(X))=(2, 4)$.

Assume that $X$ is a canonically polarized Gorenstein $3$-fold with at most terminal singularities
and with $(K_X^3,p_g(X))=(2, 4)$.
If $\dim \mathrm{Im}\phi_{K_X}=3$, then $\phi_{K_X} \colon X \rightarrow \mathbb{P}^3$ is a finite double cover of $\mathbb{P}^3$ branched along
a surface of degree $10$, according to the results of Fujita \cite{fujita1} and \cite{fujita2}.

If $\dim \mathrm{Im}\phi_{K_X}=2$, then \lemref{lem:fibration}, \lemref{lem:section} and \propref{prop:equalities} except the last assertion (f) still hold.
However, $\pi(\Gamma)$ is a point according to \propref{prop:equalities}~(d) and we won't be able to use the same techniques in the proof of \lemref{lem:alpha}.
\end{rem}

\subsection{Proof of Theorem~1.2}
We are able to complete the proof of \thmref{thm:baselocus}.
Let $X$ be a $3$-fold as in the assumption of \thmref{thm:baselocus}.
We stick to the same notation as above.

We have seen that $|K_X|$ has no fixed part, i.e., $\overline{Z}=0$ by \propref{prop:equalities}~(a).
Note that $\mathrm{Bs}|K_X|\not=\emptyset$ because $|K_X|$ is ample and $\dim \mathrm{Im}\phi_{K_X}=2$.
Also recall that from \lemref{prop:equalities}~(f) and \lemref{lem:section} that $\Gamma \le E_l|_S$
and
\begin{align}\label{eq:EVreduced}
E_V=(\sum_{i=1,i\not=l}^na_iE_i)|_S+(E_l|_S-\Gamma),\ E_V'=(\sum_{i=1,i\not=l}^{n}b_iE_i)|_S+(E_l|_S-\Gamma)
\end{align}

\begin{lem}\label{lem:basecurve}
The base locus of $|K_X|$ consists of a unique irreducible curve $\overline{\Gamma}:=\pi(\Gamma)$.
\end{lem}
\proof
Note that $\pi(\Gamma)$ is an irreducible curve and $\pi(\Gamma)=\pi(E_l)$ by \propref{prop:equalities}~(d)-(f).
Let $\overline{\Gamma}:=\pi(\Gamma)$.
Then $\overline{\Gamma} \subseteq \mathrm{Bs}|K_X|$.

Assume that $\overline{\Lambda}$ is an irreducible curve such that $\overline{\Lambda} \subseteq \mathrm{Bs}|K_X|$
and $\overline{\Lambda} \not= \overline{\Gamma}$.
Recall the construction of $\pi=\alpha \circ \beta$ from \propref{prop:fgresolution} and \propref{prop:blowupresolution}
and the fact that $\alpha(D_t)$ is a point.
We see that there exists an irreducible curve $\Lambda$ contained in some $E_i|_S$  such that $\overline{\Lambda}=\pi(\Lambda)$.
Since $\overline{\Lambda}\neq \overline{\Gamma}$, we conclude that $i\not=l$ and
$\Lambda$ is contained in $\mathrm{Supp}E_V$.
Hence $K_X \overline{\Lambda}=0$ by \propref{prop:equalities}~(e).
Therefore $\overline{\Lambda}=0$ since $K_X$ is ample, a contradiction.

Suppose $p$ is an isolated point of $\mathrm{Bs}|K_X|$ and $p \not \in \overline{\Gamma}$.
Recall that $\pi^*K_X|_S=M|_S+\Gamma+D_V'+E_V'$ by \eqref{eq:fulladjunction} and \eqref{eq:vertical}.
Therefore $\pi^{-1}(p) \cap S \subseteq \mathrm{Supp}(D_V'+E_V')$.
We may write $D_V'+E_V'=A+B$,
where $A$ and $B$ are effective divisors such that $\pi(\mathrm{Supp}A)=p$ and $p\not\in \pi(\mathrm{Supp}(B))$.
Then $A\cdot B=0$ and $A.\Gamma=0$ since $p \not \in \overline{\Gamma}$.
Because $D_V'+E_V'$ is contained in the fibers of $f|_S$ by \lemref{lem:section}, we have $M|_S.A=0$.
We also conclude from \propref{prop:equalities}~(e) that $\pi^*K_X\cdot A=0$.
Then $A^2=A(\pi^*K_X|_S-M|_S-\Gamma-B)=0$.
But $\pi^{*}K_X|_S$ is nef and big, this contradicts the algebraic index theorem.
This completes the proof.\qed.

\begin{lem}\label{lem:alpha}
The birational morphism $\alpha \colon X_0 \rightarrow X$ in \propref{prop:fgresolution} is indeed the identity morphism
and $\overline{\Gamma}$ is contained in the smooth locus of $X$.
\end{lem}
\proof
Because $\overline{\Gamma}=\pi(\Gamma)$ is a curve  and $\alpha(D_t)$ is a point for $1 \le t \le m$,
we conclude that $\Gamma$ is not contained in $\beta^*D_t$ and then,
by \propref{prop:equalities}~(b) and \lemref{lem:section}, $\Gamma\cap\beta^*D_t=\emptyset$ for any $t$.
This yields $\beta(\Gamma)\cap D_t=\emptyset$ for any $t$.
Because $\cup_{t=1}^mD_t$ is a disjoint union of fibers of $\alpha$ by \propref{prop:fgresolution}~(d) and $\pi=\alpha\circ \beta$,
$\alpha(D_t) \not \in \pi(\Gamma)$ for any $t$.
Since $\overline{\Gamma}=\mathrm{Bs}|K_X|=\mathrm{Bs}|\overline{M}|$ by \lemref{lem:basecurve},
we see that $\alpha(D_t) \not \in \mathrm{Bs}|\overline{M}|$ for any $t$.
According to \propref{prop:fgresolution}~(b), we conclude that $\alpha \colon X_0 \rightarrow X$ is the identity morphism
and $\overline{\Gamma}=\mathrm{Bs}|\overline{M}|$ is contained in the smooth locus of $X$.
\qed.

It follows from the previous lemma that $\pi=\beta$ (see \propref{prop:blowupresolution})
and \eqref{eq:fulladjunction} becomes
\begin{align}\label{eq:reducedadjunction}
K_Y=\pi^*K_X+\sum_{i=1}^na_iE_i=\pi^*K_X+\Gamma+E_V,\ \pi^*K_X=M+\sum_{i=1}^nb_iE_i=M+\Gamma+E_V'
\end{align}

For a general member $S \in |M|$, denote by $\overline{S}:=\pi(S)$ and by $\sigma:=\pi|_S \colon S \rightarrow \overline{S}$.
Note that $\overline{S}\in |K_X|$.

\begin{lem}\label{lem:Sbar}
We may choose a general $S \in |M|$ such that both $S$ and $\overline{S}$ are smooth. Moreover,
\begin{enumerate}[\upshape (a)]
    \item the exceptional divisors of $\sigma \colon S \rightarrow \overline{S}$ is contained in $E_V$;
    \item the formula $K_X|_{\overline{S}}=(p_g(X)-2)\overline{C}+\sigma_\ast(\Gamma)$ holds, where $|\overline{C}|$ is a base-point-free pencil of curves
          induced by the fibration $f|_S$, i.e., $\sigma^*\overline{C}=C$.
\end{enumerate}
\end{lem}

\proof
According to Bertini's theorem,
we may choose $S \in |M|$ such that
\begin{enumerate}[\upshape (1)]
\item $\overline{S}$ is smooth outside $\overline{\Gamma}$ by \lemref{lem:basecurve};
\item $\overline{S}$ is smooth at a general point of $\overline{\Gamma}$ by \propref{prop:equalities}~(f)
      and \propref{prop:blowupresolution}~(c);
\item $\overline{S}$ is smooth at the point $\pi_0 \circ \cdot \circ \cdot \pi_{k-1}(W_k)$ for those $k$
       such that $a_k=2b_k$ by \propref{prop:blowupresolution}~(d).
\end{enumerate}
In particular, $\overline{S}$ is a normal surface.

Note that $\sigma=\pi|_S $ is isomorphic at the points outside $\Gamma \cup \mathrm{Supp}E_V$.
To show the smoothness of $\overline{S}$, it suffices to show that $\overline{S}$ is smooth at the points where $\sigma^{-1}$ is not defined.
Let $q$ be such a point.
Then $q \in \overline{\Gamma}$ by the choice of $S$.
According to Zariski's main theorem, $\sigma^{-1}(q)$ is a connected curve.
Since $q \in \overline{\Gamma}$, there is an irreducible component $Q$ of $\sigma^{-1}(q)$
such that $Q\Gamma>0$.
Note that $Q \le E_V$ and recall that $(2E_V'-E_V)\Gamma=0$ by \propref{prop:equalities}~(b).
It follows from \eqref{eq:EVreduced} that $Q \le E_k|_S$ for some $k \not =l$ and $a_k=2b_k$.
Hence $\overline{S}$ is smooth at $q$ by the choice of $S$.

Because both $S$ and $\overline{S}$ are smooth, $\sigma$ is a composition of blowups
and the support of the exceptional divisors of $\sigma$ coincides with $\mathrm{Supp} E_V$.
Recall that $f|_S$ is the fibration induced by $|M||_S$,
and $C$ is a general fiber of $f|_S$,
and $E_V$ is contained in the fibers of $f|_S$.
Hence there is a base-point-free pencil $|\overline{C}|$ of curves such that $\sigma^*(\overline{C})=C$.

Note that $\dim |K_X||_{\overline{S}}=p_g(X)-2$.
Recall that  $\Sigma$ is a normal surface of minimal degree, i.e., $d_\Sigma=p_g(X)-2$.
Note that the target space of the fibration $f|_S$ is a hyperplane section of $\Sigma$,
which is a smooth rational curve.
Therefore $M|_S$ is linearly equivalent to $d_\Sigma C$.
Also $\sigma^*(K_X|_{\overline{S}})=\pi^*K_X|_S=M|_S+\Gamma+E_V'=d_{\Sigma}C+\Gamma+E_V'$ by \eqref{eq:reducedadjunction}
Since $\mathrm{Supp}E_V'=\mathrm{Supp}E_V$, $\sigma_\ast E_V'=0$.
We obtain $K_X|_{\overline{S}}=(p_g(X)-2)\overline{C}+\sigma_\ast(\Gamma)$.\qed

\proof[Proof of \thmref{thm:baselocus}]
We choose $S$ and $\overline{S}$ as \lemref{lem:Sbar}.

For (a), by \lemref{lem:basecurve} and \propref{prop:equalities}~(d), it remains to show that $\overline{\Gamma}$ is smooth.
Let $d$ be the degree of $\sigma|_\Gamma \colon \Gamma \rightarrow \overline{\Gamma}$.
Then $\sigma_\ast(\Gamma)=d\overline{\Gamma}$.
By \lemref{lem:Sbar}~(b), the projection formula yields $d\overline{C}.\overline{\Gamma}=\overline{C}\sigma_\ast(\Gamma)=\sigma^*\overline{C}\Gamma=C\Gamma=1$.
Hence $d=1$ and $\overline{C}\overline{\Gamma}=1$.
Because $|\overline{C}|$ is base point free, $\overline{\Gamma}$ is smooth.

For (b), let $\pi' \colon  X' \rightarrow X$ be the blowup of $X$ along the curve $\overline{\Gamma}$.
We have
    \begin{align*}
     K_{X'} = {\pi'}^*(K_X)+E,\ {\pi'}^*(K_X)=M'+E',
     \end{align*}
where $E'$ is the exceptional divisor of the blowup and $|M'|$ is the movable part of $|{\pi'}^*(K_{X})|$.
To prove (b), it suffices to prove that $|M'|$ is base point free

Let $S'$ be the strict transform of $\overline{S}$ under $\pi'$ and let $\sigma':=\pi'|_{S'} \colon S' \rightarrow \overline{S}$.
Because both $\overline{S}$ and $\overline{\Gamma}$ are smooth, $S'$ is also smooth,
${\pi'}^*\overline{S}=S'+E$ and $\sigma'$ is an isomorphism.
Moreover, $E'|_{S'}$ is a smooth rational curve and $E'|_{S'}={\sigma'}^*(\overline{\Gamma})$.

It suffices to show that the trace of $|M'|$ on $S'$ is base point free.
Note that
$M'|_{S'}={\pi'}^*(K_X)|_{S'}-E'|_{S'}={\sigma'}^*(K_X|_{\overline{S}})-{\sigma'}^*(\overline{\Gamma})$.
We have seen $\sigma_\ast(\Gamma)=\overline{\Gamma}$ and $K_X|_{\overline{S}}=(p_g(X)-2)\overline{C}+\overline{\Gamma}$ in \lemref{lem:Sbar}.
Therefore $M'|_{S'}={\sigma'}^*((p_g(X)-2)\overline{C})$.
Because $|\overline{C}|$ is base point free and $\dim|M'||_{S'}=p_g(X)-2$,
we conclude that $|M'||_{S'}$ composed with the pencil $|{\sigma'}^*\overline{C}|$.
In paticular, $|M'||_{S'}$ is base point free and complete the proof of (b).

Assertions (c) and (d) follow by \lemref{lem:fibration}~(b)-(c) and \propref{prop:equalities}~(c).\qed
\section{Classification}
The whole section is devoted to prove \thmref{thm:classification}.
We stick to the same notation in \thmref{thm:baselocus} and assume that $X$ is locally factorial (see \remref{rem:locallyfactorial}).
In particular, the left triangle of the following diagram
\begin{align}\label{eq:diagram}
\xymatrix{
Y  \ar_{\pi}"2,1" \ar^{\phi}"2,2"  \ar@{-->}^{\overline{\phi}}"1,2"     & \Sigma_e \ar^{r}"2,2"            &          \\
X  \ar@{-->}_{\phi_{K_X}}"2,2"          & \Sigma   \ar@{^{(}->}[r]    & \mathbb{P}^{p_g(X)-1}
}
\end{align}
is commutative and $Y$ is also locally factorial.
And we have
\begin{align}\label{eq:adjunction}
      K_Y = \pi^*K_X+E_0,\ \pi^*K_X=M+E_0
\end{align}
where $|M|$ is base point free and $\phi$ is induced by $|M|$.
Also $\phi|_{E_0} \colon E_0 \rightarrow \Sigma$ is a birational morphism.

According to \thmref{thm:baselocus}~(c), $\deg \Sigma=p_g(X)-2 \ge 5$.
So $\Sigma$ is obtained from a Hirzebruch surface $\Sigma_e$ for some $e \ge 0$
via the birational morphism $r \colon \Sigma_e \rightarrow \Sigma$ induced by the linear system $|s+(e+k)l|$,
where $k$ is a nonnegative integer such that
\begin{align}
p_g(X)=2k+e+2\ \text{and}\ \deg\Sigma=2k+e \label{eq:pg}
\end{align}
(see Section~2 for the notation for the Hirzebruch surfaces).
More precisely, we have two possibilities as follows
(cf.~\cite[Exercises~IV.18~4)]{algebraicsurface} or \cite[p.~380, Corollary~2.19]{gtm52}).
\begin{enumerate}[\upshape (1)]
 \item If $k \ge 1$, then $r$ is an isomorphism and $\Sigma$ is indeed smooth.
       In this case,
       $\phi|_{E_0} \colon E_0 \rightarrow \Sigma$ is indeed an isomorphism since $E_0$ is also a Hirzebruch surface.
       If $e>0$, then the ruling $|l|$ of $\Sigma_e$ coincides with
       $\pi|_{E_0} \colon E_0 \rightarrow \overline{\Gamma}$ via $r^{-1}\phi|_{E_0}$
       because $\Sigma_e$ has a unique ruling.
       We may assume this still holds when $e=0$ by possibly exchanging the two rulings of $\Sigma_0$.

 \item If $k=0$, then $e \ge 5$ and $\Sigma$ is a cone over a rational normal curve.
       Moreover, $r$ contracts exactly the negative section $s$ and $\mathrm{v}:=r(s)$ is the vertex of the cone $\Sigma$.
\end{enumerate}
To prove \thmref{thm:classification}, we shall exclude the case (2).
The following lemma allows us to treat both cases in a unified way.
\begin{lem}\label{lem:factorthrough}
The rational map $\overline{\phi}=r^{-1}\phi$ is indeed a morphism and
$\overline{\phi}|_{E_0} \colon E_0 \rightarrow \Sigma_e$ is an isomorphism.
\end{lem}
\proof
The lemma is nontrivial only for $k=0$.
In this case, since $\rho(E_0)=2$ and $\rho(\Sigma)=1$,
the birational morphism $\phi|_{E_0} \colon E_0 \rightarrow \Sigma$ contracts exactly the negative section
$s_{E_0}$ of the ruling $\pi|_{E_0} \colon E_0 \rightarrow \overline{\Gamma}$ and
maps any fiber $l_{E_0}$ of this ruling to a line on $\Sigma$.
Since $\phi|_{E_0}$ is induced by $|M||_{E_0}$, we conclude that $M|_{E_0}=s_{E_0}+el_{E_0}$.

Assume by contradiction that $r^{-1}\phi$ is not a morphism.
Then the locus where $r^{-1}\phi$ is not defined is contained in $\phi^{-1}(\mathrm{v})$.
Let $\widehat{\pi} \colon \widehat{Y} \rightarrow Y$ be the resolution of the indeterminacy of $r^{-1}\phi$
and let $\widehat{\phi} \colon \widehat{Y} \rightarrow \Sigma_e$ be the induced morphism
such that $\phi\widehat{\pi}=r\widehat{\phi}$:
\begin{align*}
\xymatrix{
    \widehat{Y} \ar^{\widehat{\phi}}"1,2" \ar_{\widehat{\pi}}"2,1" & \Sigma_e \ar^{r}"2,2" &  s=r^{-1}(\mathrm{v}) \ar"2,3"  \\
    Y           \ar^{\phi}"2,2"                                    & \Sigma                &  \mathrm{v}
    }
\end{align*}
We may assume that $\widehat{\pi}|_{\widehat{\pi}^{-1}(Y\setminus\phi^{-1}(\mathrm{v}))}\colon \widehat{\pi}^{-1}(Y\setminus\phi^{-1}(\mathrm{v}))\rightarrow Y\setminus\phi^{-1}(\mathrm{v})$ is an isomorphism and that $\widehat{\pi}^{-1}(\phi^{-1}(\mathrm{v}))$ is contained in the smooth locus of $\widehat{Y}$.

From the commutative diagram,
we have
$$\widehat{\pi}^*M=\widehat{\pi}^*\phi^*\O_{\Sigma}(1)=\widehat{\phi}^*r^*\O_{\Sigma}(1)=
\widehat{\phi}^*(s+el)=\widehat{\Delta}+e\widehat{L},$$
where $\widehat{\Delta}:=\widehat{\phi}^*s$, $\widehat{L}:=\widehat{\phi}^*l$ and $|\widehat{L}|$
is a base-point-free pencil of divisors.
It follows that $M =\Delta+eL$,
where $\Delta:=\widehat{\pi}_{\ast}\widehat{\Delta}$, $L:=\widehat{\pi}_{\ast}\widehat{L}$
and $|L|$ has no fixed part.

Recall that $C$ is a general fiber of $\phi$ and $E_0C=1$.
Since $MC=0$, we have $\Delta C=LC=0$.
Then $\Delta \not \ge E_0$ and $L \not \ge E_0$ since $E_0C=1$,

We now show that $|L|$ is a base-point-free pencil.
Note that $e+2=h^0(Y, \O_Y(M)) \ge h^0(Y, \O_Y(eL))) \ge eh^0(Y, \O_Y(L))-e+1$.
Therefore $\dim |L|=1$.
Since $\Delta|_{E_0}+eL|_{E_0}=M|_{E_0}\sim s_{E_0}+el_{E_0}$ and $\Delta \not \ge E_0$,
we have $\Delta|_{E_0} > 0$ and thus $\Delta > 0$.
Moreover, from the commutativity of the diagram above and the definition of $\Delta$,
we see that $\phi(\mathrm{Supp}\Delta)=\mathrm{v}$.
Because $\phi|_{E_0}$  contracts exactly the curve $s_{E_0}$,
$\Delta|_{E_0}=bs_{E_0}$ for some integer $b \ge 1$.
Then $e^2\pi^*K_X.L^2=(M+E_0)(M-\Delta)^2=M|_{E_0}^2-2M|_{E_0}\Delta|_{E_0}+\Delta|_{E_0}^2=e(1-b^2)$.
Since $K_X$ is ample and $|L|$ has no fixed part, we obtain $\pi^*K_X.L^2=0$, $b=1$ and
thus $L|_{E_0}=\frac 1e(M-\Delta)|_{E_0}=l_{E_0}$.
Then the trace of the pencil $|L|$ on $E_0$ is $|l_{E_0}|$ since $L \not \ge E_0$.
This implies   $\mathrm{Bs}|L|\cap E_0=\emptyset$.
If $\mathrm{Bs}|L|\not=\emptyset$, $L^2$ is rationally equivalent to an effective $1$-cycle whose support is not contained in $E_0$ since  $\dim |L|=1$.
Hence the ampleness of $K_X$ and $\pi^*K_X.L^2=0$ imply that $|L|$ is base point free.

Because both $|L|$ and $\widehat{L}$ are base point free and $L=\widehat{\pi}_{\ast}\widehat{L}$,
we have $\widehat{\pi}^*L=\widehat{L}$.

Let $F$ be any irreducible and reduced  curve contracted by $\widehat{\pi}$.
Then $\widehat{L}.F=\widehat{\pi}^*L.F=0$.
On one hand, since $\widehat{L}=\widehat{\phi}^*(l)$,
$\widehat{\phi}(F)$ is contained in one of the fiber of ruling induced by $|l|$.
On the other hand, $r^{-1}\phi$ is defined outside $\phi^{-1}(\mathrm{v})$,
so $r\widehat{\phi}(F)=\phi\widehat{\pi}(F)=\mathrm{v}$ and thus
$\widehat{\phi}(F)$ is contained in $r^{-1}(\mathrm{v})=s$.
Therefore $\widehat{\phi}(F)$ is a point in $\Sigma_e$.
This means $\widehat{\phi}$ factors though $Y$.

Hence $\overline{\phi}$ is a morphism.
Its restriction $\overline{\phi}|_{E_0} \colon E_0 \rightarrow \Sigma_e$
is birational because so is $\phi|_{E_0}$.
Then it is an isomorphism because both $E_0$ and $\Sigma_e$ are Hirzebruch surfaces.\qed

\begin{rem}\label{rem:locallyfactorial}
We  need the assumption that $X$ is locally factorial to apply intersection theory.
For example, if this assumption is dropped, then we do not know whether $\Delta:=\widehat{\pi}_{\ast}\widehat{\Delta}$ is $\mathbb{Q}$-Cartier.
In this situation, $\Delta|_{E_0}$ might be not well-defined.

This assumption is also important in the proof of \lemref{lem:flatness}.
\end{rem}

From now on, we denote by $j$ the inverse of the isomorphism $\overline{\phi}|_{E_0} \colon E_0 \rightarrow \Sigma_e$.
By abuse of notation, we identify $\mathrm{Pic}(E_0)$ with $\mathrm{Pic}(\Sigma_e)=\mathbb{Z}l\oplus\mathbb{Z}s$.
Since $M=\phi^*\O_{\Sigma}(1)$ and $r^*\O_{\Sigma}(1)=s+(e+k)l$, we have
\begin{align}
M=\overline{\phi}^*(s+(e+k)l), M|_{E_0}=s+(e+k)l \label{eq:M}
\end{align}
Since $r\overline{\phi}=\phi$, we still denote by $C$ the general fiber of $\overline{\phi}$.

\begin{lem}\label{lem:flatness}
Every fiber of $\overline{\phi}$ is $1$-dimensional, reduced  and irreducible.
In particular, $\overline{\phi}$ is flat.
\end{lem}
\proof Assume that $\Phi$ is a $2$-dimensional irreducible component of a fiber of $\overline{\phi}$.
Recall that $\overline{\phi}|_{E_0}  \colon E_0 \rightarrow \Sigma_e$ is an isomorphism by \lemref{lem:factorthrough}.
It follows that $E_0 \cap \Phi=\emptyset$ or $E_0 \cap \Phi$ consists of a single point.
Since $Y$ is locally factorial, every divisor of $Y$ is Cartier, and thus $\dim E_0 \cap \Phi=1$ if $E_0 \cap \Phi \not =\emptyset$.
We conclude that $E_0 \cap \Phi=\emptyset$.
Therefore $(\pi^*K_X)^2.\Phi=(M+E_0)^2.\Phi=0$ according to \eqref{eq:adjunction}.
Since $K_X$ is ample, it follows that $\dim \pi(\Phi) \le 1$.
Then $\pi(\Phi)=\overline{\Gamma}$ and $\Phi=E_0$, a contradiction.
Hence every fiber of $\overline{\phi}$ is $1$-dimensional
and $\overline{\phi}$ is flat (cf.~\cite[p.~179, Theorem~23.1 and Corollary]{commutativering}).

Because $K_X$ is ample and $\pi^*K_X.C=1$ by \thmref{thm:baselocus}~(c),
if $\overline{\phi}$ has reducible fibers,
then some irreducible component of some reducible fiber is contained in $E_0$.
This contradicts that $E_0$ is a section of $\overline{\phi}$.
Therefore any fiber is irreducible.
Also $\pi^*K_X.C=1$ implies that any fiber is reduced.
\qed

\begin{lem}\label{lem:pushfoward}
Let $\E:=\overline{\phi}_\ast\O_Y(2E_0)$.
\begin{enumerate}[\upshape (a)]
\item Then $\E$ is a locally free sheaf of rank $2$ and the natural morphism $\overline{\phi}^*\E \rightarrow \O_Y(2E_0)$ is surjective.

\item Let $P:=\mathbb{P}_{\Sigma_e}(\E)$ and let $\tau \colon P \rightarrow \Sigma_e$ be the natural projection.
      Then the $\Sigma_e$-morphism $\psi \colon Y \rightarrow P$ associated to
      the surjective morphism $\overline{\phi}^*\E \rightarrow \O_Y(2E_0)$ is finite of degree $2$.
\item Let $E=\psi(E_0)$.
      Then $E$ is a section of $\tau$ and it is an irreducible connected component of the branch divisor of $\psi$.
\item The section $E$ of $\tau$ corresponds to the surjective morphism
      $\E=j^*\overline{\phi}^*\E \rightarrow j^*\O_Y(2E_0)$,
      whose kernel is $\O_{\Sigma_e}$.
\end{enumerate}
\end{lem}

\proof
Let $C$ be any fiber of $\overline{\phi}$. Then $p_a(C)=2$ by \thmref{thm:baselocus}~(c).
Because $\overline{\phi}$ is a flat morphism,
$C$ is Gorenstein.
Since  $C$ is reduced and irreducible by the previous lemma,
$|\o_C|$ is base point free by \cite[Theorem~3.3]{embedding}.
We have $K_C=K_Y|_C=2E_0|_C$ from \eqref{eq:adjunction}.
So $\O_Y(2E_0)|_C$ is generated by global sections and $h^0(C, \O_Y(2E_0)|_C)=p_a(C)=2$.
Then (a) follows by Grauert's Theorem.

For (b), we have the following diagram such that $\tau\psi=\overline{\phi}$.
\begin{align*}
\xymatrix{
                                             & P \ar^{\tau}"2,2"\\
Y \ar^{\psi}"1,2" \ar^{\overline{\phi}}"2,2" & \Sigma_e \ar@/^/^{j}"2,1"
}
\end{align*}
Note the the restriction $\psi|_C$ is indeed the canonical map of $C$,
which is a finite morphism of degree $2$ to the projective line $\mathbb{P}^1$.
This proves (b).

For (c), first note that the branch divisor of $\psi$ is pure of dimension one because $Y$ is normal and $P$ is smooth.
Since $Y$ has only finite many singularities and $Y$ is locally factorial,
we conclude that the irreducible components and the connected components of the branch divisor  coincide.
Because $E_0$ is a section of $\overline{\phi}$,  $E$ is a section of $\tau$.
Moreover, since $E_0|_C$ consists one point and $K_C=2E_0|_C$,
we conclude that $E_0|_C$ is a ramification point of the canonical morphism of $C$.
Hence $E_0$ is an irreducible component and thus a connected component of the branch divisor.

Note that $E=\psi(E_0)=\psi j(\Sigma_e)$.
From the construction of $\psi$,
$\psi j$ corresponds the pullback of the surjective morphism $\overline{\phi}^*\E \rightarrow \O_Y(2E_0)$ by $j^*$,
which is $\E=j^*\overline{\phi}^*\E \rightarrow j^*\O_Y(2E_0)$.
Denote by $\mathcal{K}$ its kernel.
Then $\O_P(E)\otimes\tau^*\mathcal{K}=\O_P(1)$ (see the proof of \cite[p.~371, Propostion~2.6]{gtm52}).
Applying $(\psi j)^*$ to this equality,
since $\psi^*E=2E_0$ by (c), $\tau (\psi  j)=\mathrm{id}_{\Sigma_e}$
and $(\psi j)^*\O_P(1)=j^*\O_Y(2E_0)$,
we conclude that $\mathcal{K}=\O_{\Sigma_e}$.
\qed\smallskip

Let $D$ be the branch divisor of the double cover $\psi \colon Y \rightarrow P$.
Then
\begin{align}
D \sim 2\L\ \text{and}\ K_Y=\psi^*(K_P+\L)\label{eq:classifycoverdata}
\end{align}
for some $\L \in \mathrm{Pic}(P)$.
Since $P=\mathbb{P}_{\Sigma_e}(\E)$ with $\E$ as an extension
\begin{align}\label{eq:classifyrank2}
0 \rightarrow \O_{\Sigma_e} \rightarrow \E \rightarrow j^*\O_Y(2E_0) \rightarrow 0
\end{align}
we have $\mathrm{Pic}(P)=\mathbb{Z}E\oplus \tau^*\mathrm{Pic}(\Sigma_e)$,
\begin{align}
\O_P(E)=\O_P(1)\ \text{and}\ K_P=\tau^*(K_{\Sigma_e}+j^*\O_Y(2E_0))-2E \label{eq:classifypicard}
\end{align}
We shall determine $j^*\O_Y(2E_0)$ and $\L$ in terms of $\mathrm{Pic}(\Sigma_e)$ and $\O_P(E)$.

\begin{lem}\label{lem:L}
We have
\begin{align}
j^*\O_Y(2E_0)\cong \O_{\Sigma_e}(-2s-2al)\ \text{and}\ \L=3E+\tau^*(5s+5al), \label{eq:L}
\end{align}
where $a$ is an integer such that
\begin{align}
k=3a-2e-2         \label{eq:ka}
\end{align}
\end{lem}
\proof
By \eqref{eq:adjunction}, \eqref{eq:M} and the adjunction formula,
we have $\O_{E_0}(E_0)=\frac 13 (K_{E_0}-M|_{E_0})=-s-al$ with integer $a=\frac 13 (k+2e+2)$.
Since $j \colon \Sigma_e \rightarrow Y$ factors through $E_0$,
$j^*\O_Y(2E_0) \cong \O_{\Sigma_e}(-2s-2al)$.

Since $\psi^*E=2E_0$, according to \eqref{eq:adjunction} and \eqref{eq:M},
$K_Y = \psi^*(\tau^*(s+(e+k)l)+E)$.
On the other hand,
by \eqref{eq:classifycoverdata} and \eqref{eq:classifypicard},
$K_Y=\psi^*(\tau^*(-4s-(e+2a+2)l)-2E+\L)$.
It follows that $\psi^*(\L_0)=\O_Y$, where $\L_0=\L-(3E+\tau^*(5s+5al))$.
Note that $\psi_{\ast}\O_Y=\O_P \oplus \L^{\vee}$.
The projection formula yields $\L_0 \oplus (\L_0 \otimes \L^{\vee})=\O_P \oplus \L^{\vee}$.
It is clear that $H^0(P, \L^{\vee})=0$ and $H^0(P, \L_0 \otimes \L^{\vee})=0$.
We obtain $\L_0=\O_P$ and the required formula for $\L$.
\qed

\begin{prop}\label{prop:semipositivity}Let $\omega_{Y/\Sigma_e}=K_Y-\overline{\phi}^*K_{\Sigma_e}$.
\begin{enumerate}[\upshape (a)]
\item Then $\overline{\phi}_\ast\omega_{Y/\Sigma_e}=\E \otimes \O_{\Sigma_e}(3s+3al)$.
\item The pair $(e, a)$ with $a$ defined by \eqref{eq:ka} satisfies $a \ge e \ge 3$; or $1 \le e \le 2$, $a \ge e+1$; or $e=0$, $a \ge 2$.
\end{enumerate}
\end{prop}
\proof
Let $\o_{Y/\Sigma_e}=K_Y-\overline{\phi}^*K_{\Sigma_e}$
and let $C$ be any fiber of $\overline{\phi}$.
Then $\o_{Y/\Sigma_e}|_C=2E_0|_C=K_C$ by \eqref{eq:adjunction}.
Therefore $\overline{\phi}_{\ast}\o_{Y/\Sigma_e}$ is a locally free sheaf of rank $2$
by Grauert's theorem.
It follows that $\overline{\phi}_{\ast}\o_{Y/\Sigma_e}$ is semi-positive by \cite[Theorem~III and (1.3)~Remark~(iii)]{viehweg}.

We have seen that $K_Y=\psi^*(\tau^*(s+(3a-e-2)l)+E)$ and $\O_P(E)=\O_P(1)$.
Also $\psi_{\ast}\O_Y=\O_P \oplus \L^{\vee}$.
Applying the projection formula to $\psi$ and then to $\tau$, we obtain
$\overline{\phi}_{\ast}(\o_{Y/\Sigma_e})=\E \otimes \O_{\Sigma_e}(3s+3al)$.
It has $\O_{\Sigma_e}(s+al)$ as a quotient by \eqref{eq:classifyrank2} and \eqref{eq:L}.
So $a-e=\deg \O_{\Sigma_e}(s+al)|_s \ge 0$ by the semi-positivity of $\overline{\phi}_{\ast}\o_{Y/\Sigma_e}$.

Note that $p_g(X)=6a-3e-2$ by \eqref{eq:pg} and \eqref{eq:ka}.
Since $p_g(X) \ge 7$ by the assumption of \thmref{thm:classification},
it is clear that the pair $(e, a)$ satisfies the required inequalities.\qed

\proof[Proof of \thmref{thm:classification}]
By \propref{prop:semipositivity}, we see that $k=3a-2e-2 \ge 1$ and thus $\Sigma$ is the embedding of $\Sigma_e$ in
$\mathbb{P}^{p_g(X)-1}$ by the discussion at the beginning of this section.
Moreover, we have seen $p_g(X)=6a-3e-2$ in the proof of the proposition above and thus $K_X^3=8a-6e-4$
by the Noether equality. Hence (a) is established.

Assertion (b) follows from \lemref{lem:factorthrough} and \lemref{lem:flatness} since $r\circ \overline{\phi}=\phi$.

Because $j^*\O_Y(2E_0)\cong \O_{\Sigma_e}(-2s-2al)$ by \lemref{lem:L} and $a \ge e$ by \propref{prop:semipositivity},
the exact sequence \eqref{eq:classifyrank2} splits by \lemref{lem:rank2}~(b) and
thus $\mathcal{E}=\O_{\Sigma_e}\oplus \O_{\Sigma_e}(-2s-2al)$.
By \eqref{eq:classifycoverdata}, \eqref{eq:L} and \lemref{lem:pushfoward}~(c),
we see that $D=E+T$, with $T \in |5E+\tau^*(10s+10al)|$
and $E\cap T=\emptyset$.
We conclude (c) and (d) directly from the discussion above, by identifying  $\overline{\phi}$ with $\phi$ via the isomorphism $r$.

If $X$ is smooth, so is $Y$. Therefore the branch locus of $\psi$ is smooth.
Comparing \eqref{eq:classifycoverdata}, \eqref{eq:classifyrank2} and \eqref{eq:L} with
\eqref{eq:rank2}-\eqref{eq:coverdata},
we conclude that $X$ is one of the $3$-folds constructed in Section~2 and complete the proof of \thmref{thm:classification}.\qed

\paragraph{Acknowledgement.}
Both authors would like to thank Meng Chen for many valuable suggestions and helpful discussions.
The first author thanks Xiaotao~Sun for all the support during the preparation of the paper.
We thank the anonymous referee for his/her valuable suggestions.
The first author was partially supported by the China Postdoctoral Science Foundation (Grant No.:~2013M541062)
and the second author is partially supported by the National Natural Science Foundation of China (Grant No.:~11171068).

\medskip

\noindent\textbf{Authors' Addresses:}\\\smallskip

\noindent Yifan~Chen,\\
School of Mathematics and Systems Science, Beijing University of Aeronautics and Astronautics, \\
Xueyuan Road No.~37, Beijing 100191, P.~R.~China\\
Email:~chenyifan1984@gmail.com\\\smallskip

\noindent Yong~Hu,\\
School of Mathematical Sciences, Fudan University,\\
Shanghai 200433, P.~R.~China\\
Email:~11110180002@fudan.edu.cn
\end{document}